\documentclass[11pt]{amsart}
\usepackage{amsmath,amsthm,amssymb,amscd, graphicx, wrapfig}
\usepackage[all]{xy}
\usepackage[margin=1.1in]{geometry}
\usepackage[usenames, dvipsnames]{color}
\usepackage{comment}
\usepackage[bookmarks=false,colorlinks,linkcolor=blue,citecolor=blue]{hyperref}

\numberwithin{equation}{section}

\newtheorem{theorem}{Theorem}

\newtheorem{exercise}[theorem]{Exercise}

\newtheorem{defn}[theorem]{Definition}
\theoremstyle{definition}
\newtheorem{exmp}[theorem]{Example}

\newtheorem{prob}[theorem]{Open Problem}
\newtheorem{exer}[theorem]{Exercise}

\title[Planes and Axioms ]{Counting plane arrangements via oriented matroids}

\author{Stefan Forcey} \address[S. Forcey]{
    Department of Mathematics\\
    The University of Akron\\
    Akron, OH 44325-4002
    }
    \email{sforcey@uakron.edu}  \urladdr{http://www.math.uakron.edu/\~{}sf34/}
    
\begin{document}

\begin{abstract}
Planes are familiar mathematical objects which lie at the subtle boundary between continuous geometry and discrete combinatorics. A plane is geometrical, certainly, but the ways that two planes can interact break cleanly into discrete sets: the planes can intersect or not. Here we review how oriented matroids can be used to try to capture the combinatorial aspect, giving a way to encode with finite sets all the ways that $n$ planes can interact. We mention how the one-to-one correspondence breaks down in 2 dimensions for 9 lines, and in 3D for 8 planes. We include  illustrations of all the types of plane arrangements using $n=4$ and 5. 
\end{abstract}


\baselineskip=17pt
\maketitle

\section{Introduction: Linear algebra}

The number of solutions to a system of linear equations can only be 0, 1 or $\infty.$ The geometric explanation of this fact is that each linear equation in $d$ variables (with a constant term) determines an affine hyperplane in $\mathbb{R}^d$: a point in $\mathbb{R}^1$, a line in $\mathbb{R}^2,$ a plane in $\mathbb{R}^3$, and so on\footnote{
The term ``hyperplane'' by itself in a book on linear algebra often means a subspace of dimension $d-1$ in $\mathbb{R}^d$, therefore containing the origin.  For work that focuses on affine geometry, as in the remainder of this paper, hyperplanes can be assumed to be affine
unless the context indicates otherwise. Specifically ``lines'' in $\mathbb{R}^2,$ and ``planes'' in $\mathbb{R}^3$  refer to either subspaces or  affine spaces.
}.
For instance, several distinct planes in  $\mathbb{R}^3$ can: (1) avoid mutual intersection, (2) mutually intersect in one point (we'll need at least three for that), or (3) have a $1$-dimensional mutual intersection. 

\begin{figure}[h]
    \centering
    \includegraphics[width=\textwidth]{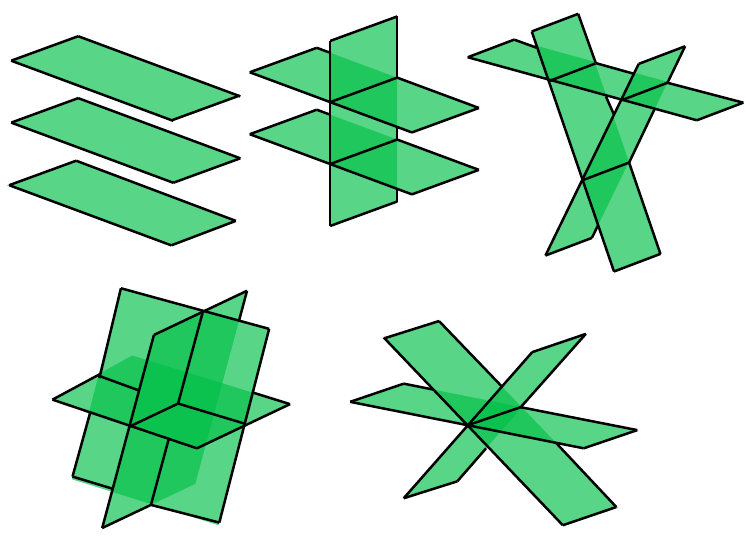}
    \caption{Here are the five ways that three planes can intersect. The top row shows the non-central cases in which there is no simultaneous solution to the three (affine) linear equations. The top left is the trivial arrangement (rank 1), the bottom left is essential (full rank, $r=3$), and the other three are rank 2.}
    \label{smallfive}
\end{figure}

A collection of $n$ distinct hyperplanes in $\mathbb{R}^d$ is called an  \textit{arrangement} of hyperplanes. If we need to stress that they include affine hyperplanes, we call it an affine arrangement. Plane arrangements and line arrangements can refer to either type, depending on context.
Recall that a system of equations which includes the origin as a solution is called \textit{homogeneous}. The corresponding \textit{central} arrangements have a nonempty intersection of all the hyperplanes which includes the origin. The \textit{essential} arrangements are the ones which have only one common point of intersection, at the origin. These essential arrangements thus correspond to \textit{full rank} homogeneous systems where $r=d$: the number of linearly independent columns of coefficients equals the number of variables.   \textit{Affine central}  arrangements are allowed to have that nonempty intersection away from the origin. \textit{Affine essential} arrangements are full rank, and thus have some subsets of their hyperplanes possessing a single point of intersection. Note that non-central always implies affine. Once again, when working in the affine world, we often drop the adjective and refer simply to essential and central arrangements.  The distinction of affine, or not, is an important feature in linear algebra, since systems are easier to solve when they are homogeneous,  and since an affine system which has a solution can be solved in stages via a detour through the associated homogeneous version. Similarly, we will see how any \textit{non-trivial} affine arrangement can be lifted to an essential one called its suspension, or cone. 
Images of two parallel or intersecting lines and three parallel or intersecting planes---
which have 5 different possibilities, shown in Figure~\ref{smallfive}---are familiar from the first pages of books on linear algebra.

\subsection{Overview} Most arrangements have no common point, and so correspond to systems of linear equations with no solution. This starts with three planes: in Figure~\ref{smallfive} we see that three out of the five possible arrangements have no common intersection.
For several years it has been a mystery to this author why the pictures stop there: why are there no collections of pictures showing all the ways that four or five or more planes can interact? We fix that here: see Figures~\ref{nineof14}--\ref{fivetriple}. Making sure that our illustrated listings are correct was the motivation for this paper. That requires clearly distinguishing the arrangements, using geometric combinatorics that we review in the next section. To ensure we find them all, we need a countable model. With that goal in mind, operations on the pieces of an arrangement are described in a section on geometric algebra. Following that, we present a gentle introduction to the abstract algebraic axioms of oriented matroids.  In the section on counting we review the historical discoveries of how well these algebraic objects serve to model the combinatorial classes (from the late 1970's and early 80's). The correspondence between oriented matroids and plane arrangements breaks down in 2 dimensions for 9 lines, and in 3D for 8 planes, leaving us the consolation prize of upper bounds. Something about the geometry gets in the way of a clean representation via strings of symbols---at least as far as our current understanding can tell. We hope to motivate some improvements!

\begin{table}[ht!]
\begin{tabular}{|c|ccccccccccc|}
\hline
$r=$& $n=$& 1&2&3&4&5&6&7&8&9&10 \\
\hline
1&&1&1&1&1&1&1&1&1&1&1\\
2&&&1 & 3&8&46&790&37829&4134939&?&?\\
3&&&&1& 5& 27& 1063& 1434219&?&?&?\\
\hline
$r\le 3$&& 1& 2 & 5 &14&74&1854&1472049&?&?&?\\
  \hline
\end{tabular}  \vspace{.12 in}
\caption{ Numbers of isomorphism equivalence classes of loop-free affine oriented matroids with $n$ elements and rank $r\le 3$, also called abstract dissection types, as found in \cite{finschi}. The  values in the last row are also the numbers of face-combinatorial equivalence classes of plane arrangements in $\mathbb{R}^3$, for $n\le 7$. }\label{counts}
\end{table}

Finally, the only new content of this paper is in our last two sections where we summarize the potential ways for 4 or 5 planes to be arranged in $\mathbb{R}^3.$ We use numbers calculated by Lukas Finschi in his Ph.D. thesis \cite{finschi} to count them (see Table~\ref{counts}) and then illustrate them  in order to explicitly show that the abstract possibilities are all realized by actual planes. We have some heuristic guidance for that effort, specifically to focus on numbers of faces, from the Ph.D. thesis of Thomas Zaslavsky \cite{zas}. The success of our effort is predicted for $n\le 7$ by Goodman and Pollack \cite{proof}. Some interesting corollaries are seen, like the fact that choosing one of 10 polytopal complexes with 2 or 3 chambers determines an  arrangement of 5 hyperplanes in $\mathbb{R}^3;$ there are 64 other arrangements determined by other facts.  We include some exercises and open problems.

\subsection{Acknowledgments} The author wishes to thank the very helpful reviewers, whose assistance vastly improved the style and clarity of the exposition. Thanks also to Thomas Zaslavsky for taking the time to offer many important suggestions for correcting and polishing the penultimate draft. Any imperfections that remain are mine alone. 

\section{Combinatorics: comparing arrangements}\label{comb}

By five distinct ways of intersecting 3 planes in $\mathbb{R}^3$ we mean five arrangements of planes with \textit{non-isomorphic face posets}. We consider two arrangements of (affine) hyperplanes in $d$-dimensional Euclidean space to be equivalent if they have \textit{isomorphic face posets} (or \textit{face lattices}\footnote{All we need here is the face poset structure, but for further discussion of face lattices see \cite{obook}.}).  Here is the intuitive meaning: A collection of planes in space partition the space into regions and also partition each other via intersections. We label all those regions, sections of planes, and sections of intersections partitioned by further intersections, and call them labeled \textit{faces}. 0-dimensional faces are usually called points,
1-dimensional faces are usually called edges, and 2-dimensional faces are often called
faces or 2-cells (or a chamber if $d = 2$).  We can list those faces in a hierarchy of boundaries. In general, the highest-dimensional faces are regions of the ambient space; these are called \textit{chambers} and have boundaries made of subfaces, pieces of the hyperplanes. Those pieces are in turn bounded by lower dimensional subfaces, and so on.  The set of faces is partially ordered by  boundary subface inclusion: 
 We write $A <  B$ to mean that a
 face $A$  is included in the boundary of $B.$

\begin{figure}
     \centering
     \includegraphics[width=0.85\linewidth]{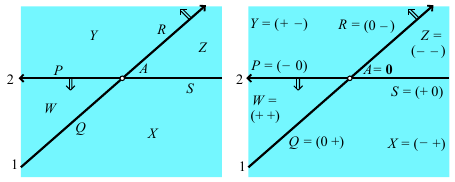}
     \caption{On the left is a line arrangement with its faces labeled: chambers $X, Y, Z, W,$ edges $P, Q, R, S,$ and a point $A.$ On the right the same faces are shown with their associated sign vectors; the small double arrows point to the plus side of each line. These will be discussed in Example~\ref{newsmole} and in the following section on oriented matroids. }
     \label{fig:newsmole}
 \end{figure}

For example, we label some faces in the line arrangement of  Figure~\ref{fig:newsmole}: the chambers $X, Z, W,$ and $Y$, edges $S, Q, R,$ and $P$, and point $A.$ 
 For faces in Figure~\ref{fig:newsmole}, $A<P<Y.$ Notice that the rank $r$ of an arrangement can be defined as the maximum number of strict inequalities in an increasing chain of related faces, so $r=2$ in Figure~\ref{fig:newsmole}. It is straightforward to see that this rank is the same as the rank of a system of linear equations describing it!  

\begin{figure}
     \centering
     \includegraphics[width=0.75\linewidth]{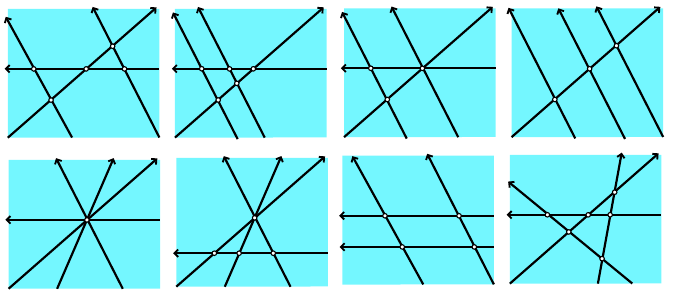}  
     \caption{Eight non-trivial affine arrangements of four lines. At bottom left is the only central arrangement. At bottom right is the only arrangement in general position.}
     \label{fig:new8}
 \end{figure}

Two posets are isomorphic when there is a bijection $f$ between their elements, and $f$ and its inverse both preserve the partial ordering. Isomorphism is an equivalence relation. Counting the number of these face-combinatorial equivalence classes of arrangements, given $n$ distinct (affine) hyperplanes, is definitely a hard open problem, even for lines. We note that rank is clearly an invariant of poset isomorphism.  In Figure~\ref{fig:new8} we show 8 nontrivial line arrangements, all of rank $r=2$, one for each face-combinatorial equivalence class.  Verifying that they are all inequivalent is straightforward: Since a poset isomorphism preserves all the ordering, any listing of the faces by dimension, boundedness, and shape will be an isomorphism invariant.
 For instance, a bounded quadrilateral chamber is a maximal element of the face poset that has 4 subfaces immediately less than it in the face ordering. This feature would be preserved by any poset isomorphism, so the first  two arrangements in Figure~\ref{fig:new8} cannot be equivalent; only the second has a bounded quadrilateral face. Unbounded chambers and edges are also faces and contain subfaces. However, an unbounded quadrilateral could not be taken to a bounded one by an isomorphism, since the number of 0-dimensional subfaces is three for the former and four for the latter.

 Finding formulas for the number of face-combinatorial equivalence classes of affine arrangements in any dimension $d\ge 2$ is an open question. Note that determining the number of 3D polyhedra with $n$ planar faces is a well-known open subproblem.
Some partial answers to the question exist. For instance, sequence A241600 in \cite{oeis} counts the number of arrangements of lines in the affine plane, up to $n=7,$ agreeing with the counts in \cite{finschi, finschiart}. Peter Shor points out that sequence A241600 is defined differently than via face-combinatorial equivalence, rather it uses parameterized equivalence, where the homotopy between equivalent arrangements must preserve straight lines and their intersections \cite{shor}.

\begin{figure}[b]
     \centering
     \includegraphics[width=\textwidth]{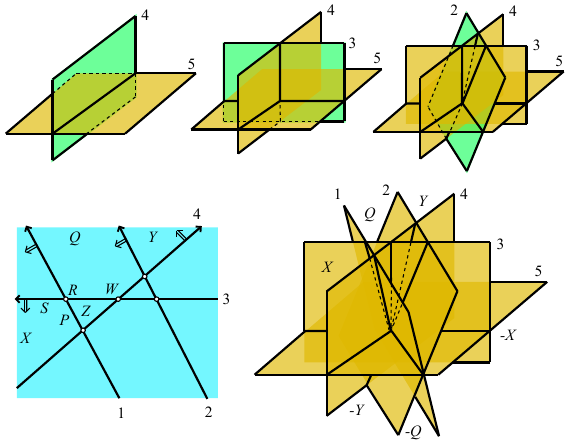}
     \caption{The top row shows the construction of a plane arrangement by adding one (highlighted) plane at a time, ending with the essential arrangement of five planes below right. The affine arrangement of four lines shown at left corresponds to the marked arrangement where the horizontal plane (5) is marked. Thus it can be seen as a top down view, or cross section on the positive side of plane 5. Chambers include $X,Q,Y,Z$; faces $S$ and $P$ are 2D but are seen as edges from above, and faces $R$ and $W$ are edges but seen as points from above.  From the line arrangement we can recover the 3D version by taking the suspension. 
     Compositions and restrictions of the labeled faces in 3D are seen in Example~\ref{examplo}.
     }
     \label{fig:fivedu}
 \end{figure}

 The lower right picture of Figure~\ref{fig:fivedu} is an example of an essential plane arrangement of 5 planes in $\mathbb{R}^3.$ It is a \textit{suspension}, or \textit{cone} of the affine line arrangement to its left. The suspension is found by placing the lower dimensional (here, 2D) arrangement in the Euclidean affine hyperplane at constant height on  the $n^{th}$ axis (here $z=1$), and then taking as the new hyperplanes the spans of the original lines, plus the plane at $z=0.$
 The suspensions of two different line arrangements  can be equivalent. This highlights the difficulty of counting the actual number of ways to arrange $n$ planes in $\mathbb{R}^3.$ 

\section{Geometric Algebra: Operations on faces of arrangements}\label{alg}

To grasp the number of ways a collection of hyperplanes might interact, mathematicians have looked hard at the  poset of faces to see if it has some additional features which correspond nicely to geometry. Here, we describe two such useful structures, with examples in Example~\ref{newsmole}, which refers to Figure~\ref{fig:newsmole}.
First there is a product structure called composition, where any two faces $X$ and $Y$ produce a third face $Z =X\circ Y.$ 
This composition can be described geometrically. The interior of a face $X$ is defined to be the points of $X$ which are not part of any proper bounding subface $Y < X$. Consider the line segment between two distinct chosen points: the first being any point in the interior of $X$ and second being any point in the interior of $Y.$ For $\epsilon > 0$ small enough, all the points distance $\epsilon$ from the first point along this line will lie in the same face, and that face is defined to be $Z=X\circ Y.$ In the case where $W$ is a point, we let $W\circ W = W.$
 Either we see $X\circ Y = X,$  or the composition is a higher dimensional face $Z=X\circ Y$ with $X<Z$.   This composition is not commutative; $W = Y \circ X$ is a different face in most cases. In general this composition fails to be a group operation because inverse faces usually fail to exist.  However, if we restrict our attention to the (affine) central arrangements, the common intersection is an identity element.

Another structure on the faces of an arrangement is similar  to the composition, but is only defined for certain ordered pairs of faces. If $X$ and $Y$ are completely separated by at least one (affine) hyperplane, but $Y$ is a part of any (affine) hyperplanes that $X$ is part of, we can define the restriction $X_Y$ of $X$ towards $Y$ as the subset of boundary faces of $X$ on the ``side towards'' $Y.$  That is, $X_Y$ is the set of boundary faces of $X$ whose points lie on a line connecting some point in the interior of $X$ with some point in the interior of $Y$.  Examples of both structures are illustrated in Example~\ref{newsmole}, referencing Figure~\ref{fig:newsmole}; and  Example~\ref{examplo}, referencing Figure~\ref{fig:fivedu}. 
 
\begin{exmp}\label{newsmole}
    In Figure~\ref{fig:newsmole} we show some faces of a small line arrangement. Composition  examples include $X\circ R = X$, $R\circ X = Z$, $P\circ Q = W$, and $Q\circ A= Q = A\circ Q.$ Note that $A$ is an identity element for composition. Restriction examples include $Z_Y= \{R\}$, $Y_X = \{P,A,R\}$, and $P_S = \{A\}.$
\end{exmp}


\section{Oriented matroids: abstract algebraic operations }\label{matroids}
 Studying the structures (composition, restriction) on the faces of an (affine) hyperplane arrangement allows us to abstract it by defining a certain kind of collection of vectors on a set, called an \textit{oriented matroid}. We will simultaneously describe the data and structure of oriented matroids while explaining how each component can be created naturally from an arrangement, with examples. We will end by defining properties of oriented matroids quite independently from geometry. We point out that every arrangement provides a model of the combinatorial definition, but that this function from arrangements to oriented matroids is not onto: the axioms are more general than that. 
 
 We start with a  finite \textit{ground set}, usually $[n] = \{1,\dots n\}$ which corresponds to the planes, arbitrarily numbered 1 to $n.$ Then we have a set of $n$-tuples, called \textit{sign vectors}, which are ordered lists of length $n$ made of the symbols $+,-,0.$ It is straightforward to create a sign vector from each face of a plane arrangement: Each plane is given an orientation arbitrarily, and each face $X$, of any dimension, can be described by a sign vector also called $X$.  The value of the component $X_i$ is determined by whether it sits on plane $i$ (in which case $X_i=0$), or on the plus side, or on the minus side. For a complete set of small examples see Figure~\ref{fig:newsmole}, where for instance edge $Q$ is on line 1, but on the positive side of line 2, so $Q=(0~+).$ For some 3D examples, see Example~\ref{examplo}.
 
 The goal is to recognize the sets of sign vectors which arise from arrangements. With that aim we define structures and properties that those structures must obey. 
 Sign vectors, like faces, are also partially ordered. We say $X\leq Y$ when $X_i\ \neq 0 \Rightarrow X_i = Y_i.$ Therefore, the strict inequality $X< Y$ means that $X$ is made by turning some of the non-zero components of $Y$ into 0's. This does agree with the ordering of faces from their inclusion in boundaries, but is more general.  We also define the opposite, or negative $-X$, of a vector $X$ in the obvious way: switch all the plus signs to minus, and vice versa, leaving 0's as is. That is, $(-X)_i = -(X_i)$ for $i=1,\dots, n.$
 The composition and restriction of any two sign vectors are defined as follows:
  \[
 (X\circ Y)_i= \Biggl\{\begin{aligned}
     X_i~~~,~~~ &  ~~~X_i \neq 0\\
     Y_i~~~,~~~ & ~~~X_i = 0
 \end{aligned}
\]

 \[
 X_Y= \{Z < X ~|~ X_i \neq -Y_i  \Rightarrow Z_i= X_i\}.
 \]

We can see that these operations are designed to mirror the geometric definitions of composition and restriction.   For simple examples, rework Example~\ref{newsmole}, using the vectors in  Figure~\ref{fig:newsmole}, where for instance $R\circ X = (0~-)\circ (-+) = (--)$. For an example of restriction:  $Z_Y=(--)_{(+-)} = \{(0~-)\}.$
\begin{exmp}\label{examplo}
     For a 3D case with $n=5$ we show some faces in Figure~\ref{fig:fivedu}. We let $X=(+++++), ~~Y=(---++),$ and $ ~~W=(-+0~0~+)$ so that all three vectors correspond to faces on the lower right of Figure~\ref{fig:fivedu}, (using the orientations shown on the lower left, interpreted as a top-down view).  All three faces are on the positive side of plane 5. For instance, $W$ is on the positive side of planes 5 and 2, the negative side of plane 1, and directly on the planes 3 and 4. Next, we can use the definitions of our operations on sign vectors: \end{exmp}

 $W\circ Y = (-+-++) = Q.$

 $Y\circ W = (---++) = Y.$

 $W\circ X = (-++++) = Z.$

 $X_W = \{(0~++++)\} = \{P\}.$

$X_Y = \{(0~++++),(++~0~++),(0~+~0~++),\\
\text{\hspace{.65 in}}(+~0~+++),(0~0~+++),(+~0~0~++),(0~0~0~++)\}\\ \text{\hspace{.4 in}} = \{P,S,R, \dots\}.$ \\
 
 Notice that usually the sign vector operations give the expected geometric result, but four of the sign vectors in $X_Y$ as defined abstractly are not actually seen as faces in the arrangement! The three that are seen, $P,S,$ and $R,$ fit our geometric description. With that in mind we  list the axioms of an oriented matroid, defined purely as a set of sign vectors without reference to geometry:

\begin{defn}
   An oriented matroid  on $E= [n]$ is a set $V$ of length $n$ sign vectors:  $V\subseteq \{+,-,0\}^E,$ obeying for all $i\in E$ and $X,Y \in V$ :
   
   (SV0) ~~$\mathbf{0} = (0,\dots,0)\in V.$
   
   (SV1) ~~$-X \in V.$

   (SV2) ~~$X\circ Y \in V.$ 

   (SV3) $((X_i = 0 \Rightarrow  Y_i = 0)$ and  $(\exists j \in E$ s.t. $X_j=-Y_j \neq 0)) \Rightarrow  X_Y\cap V \neq \emptyset.$

\end{defn}

For some basic examples, note that $V=\{\mathbf{0}\}$ and $V=\{+,-,0\}^E$ are oriented matroids for all $n.$ Thus the size of $V$ is between 1 and $3^n.$ 
The first two axioms are for convenience, requiring the oriented matroid to contain $\mathbf{0}$ and all opposite vectors. Axioms $SV2$ and $SV3$ say that compositions and restrictions must exist in $V$ (always for composition, but restriction only for certain ordered pairs).  Axiom $SV3$ we include here is equivalent to the one used by Edmonds and Mandel \cite{mandel}, and labeled $V3''$ in \cite{obook}. Other versions of axiom $SV3$ are used elsewhere, as in \cite{chapter}, included in \cite{goodman}.

An \textit{isomorphism} $f:V\to V'$ of oriented matroids is defined as a bijection on sign vectors induced by a relabeling of $E$ (also a bijection called $f$), together with a  reorientation of the sign vectors. Relabeling  means that the entries in every vector are permuted in the same way, while reorienting means that for some subset $\mathcal{O} \subseteq E$, all the signs indexed by that subset are switched in every vector. That is, the two functions obey:  \[f(X)_{f(i)} =   
 \Biggl\{\begin{aligned}
     -(X_i)~~~,~~~ &  ~~~i\in\mathcal{O} \\
     X_i~~~,~~~ & ~~~i\notin\mathcal{O}.
 \end{aligned}
\]
Notice (or check as a exercise) that such an isomorphism will automatically respect the poset structure, as well as the compositions and restrictions.

\subsection{Affine oriented matroids} We emphasize that the set of
sign vectors for a non-central arrangement of hyperplanes does not itself form an oriented matroid since it violates axioms $SV0$ and $SV1$. To recapture all the affine cases we need a variation: An \textit{affine oriented matroid} is defined as an oriented matroid $(E, V)$ with a choice of a \emph{marked} element of $g \in E.$  Geometrically, choosing an element to mark is the inverse of the  suspension operation discussed above. The marked element corresponds to the new plane that is inserted in the process of  suspension, and we can visualize the passage to a lower-dimensional affine arrangement by looking at an above view (from the positive side of the marked plane) or a cross section taken parallel to that marked plane. See Figure~\ref{fig:fivedu} for an example. Practically, this definition of affine oriented matroid is an efficient way to achieve the same result as simply dropping the first two axioms: a marked element can be used to select only the sign vectors for which its component is positive, and then ignored. The isomorphisms of affine oriented matroids must preserve that marking: Two affine oriented matroids $(E, V , g)$,
$(E', V', g')$ are called \textit{affine isomorphic} if there exists an isomorphism $f$ between the oriented matroids $(E, V)$ and 
$(E', V')$ for which $g' = f(g).$

 Oriented matroids also describe other mathematical structures, like directed graphs and zonotopes (certain polytopes with parallel facets).  However the important thing about them here is that they are precisely represented by arrangements of hyperplanes---with an important caveat. Every face-combinatorially unique central hyperplane arrangement (after labeling the planes $1,\dots,n$) yields a unique (up to isomorphism) oriented matroid via creating a sign vector from each face, as pictured in Figure~\ref{fig:newsmole}. However, there are extra oriented matroids, for which you can still draw an arrangement, but for which some of the planes must be replaced by curved surfaces! We show some well-known arrangements that serve as counterexamples in the next section. It is an open question whether there is some extra requirement of the vectors that will eliminate these extra oriented matroids, leaving only the ones that can be represented by perfectly flat planes. That problem is tempting but probably very hard: the first guess would be to look for some finite list of forbidden sub-arrangements whose presence would obstruct any possible representation via flat hyperplanes. However, this finite list has  been shown not to exist; see \cite{ard2} for more details.

 Since we don't have that answer yet, one way to use the (affine) oriented matroids for counting (affine) hyperplane arrangements is to produce an upper bound. First we exhaustively find all the isomorphism classes of affine oriented matroids for a given $n$ and dimension $d$, and then we try to reach that upper bound by producing each arrangement explicitly.

\section{Counting affine arrangements: Representability and Stretchability  }\label{count}

A \textit{loop} in an oriented matroid is an element $k\in E$ of the ground set where each sign vector in that oriented matroid has that $k^{th}$ entry equal to 0. The Topological Representation theorem of Folkman and Lawrence, \cite{folk},  says that isomorphism classes of loop-free oriented matroids are in  bijection with equivalence classes of  arrangements of \textit{pseudohyperplanes}. The latter, including pseudolines and pseudoplanes, are deformations of straight lines and planes but are required to obey the usual laws of intersection: for instance, two pseudolines can intersect at most once. 

This leads us to a pair of counting problems. Counting the number of (affine) oriented matroids on $[n]$, so the number of (affine) arrangements of $n$ pseudohyperplanes, is hard in itself. An open problem is to find a good formula.  Many of the results for small $n$ values were calculated in 2001, in the PhD thesis of Lukas Finschi \cite{finschi}.
 For instance, he showed that in $\mathbb{R}^3$ there are 14 arrangements of 4 planes and 74 arrangements of 5 planes. Finschi's counting algorithms rely on additional graph-theoretic interpretations of oriented matroids\footnote{  Finschi's thesis uses \textit{cocircuit graphs} and \textit{chirotope representations}. It also discusses algorithmic complexity: the necessary number of basic operations for each counting algorithm and the resulting run times.}.  The oriented matroids are easier to count first, and then he counts the affine versions by choosing marked elements. 
 
For instance, for counting the line arrangements of 4 lines in 2D, Finschi's program generates all the oriented matroids of 5 elements of rank 3. Then for each he  gets five affine oriented matroids, by choosing all possible elements to mark. Thus, pictorially speaking, the 8 arrangements in Figure~\ref{fig:new8} are found by choosing one of the four central essential arrangements from Figure~\ref{fivesus}, picking one of the 5 planes to mark, and then looking at a cross section parallel to that plane. Notice that there are 20 total ways to generate an affine oriented matroid just starting with the 4  central essential arrangements, so there will be potentially multiple ways to get any one of the 8 line arrangements. Thus the program sorts them into affine oriented matroid isomorphism classes.
We list some more results in Table~\ref{counts}, with values taken directly from Table 8.1 on page 165 of \cite{finschi}, but listed by number of planes and with a new row of totals. The rows are arranged by rank $r$, which is the least dimension needed to realize the affine oriented matroid as an affine arrangement.

\subsection{Counterexamples} 
\begin{figure}[h]
    \centering
    \includegraphics[width=\textwidth]{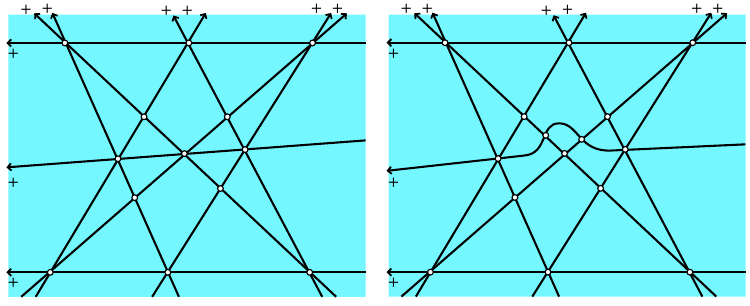}
    \caption{The Pappus arrangement on the left, and a nonstretchable pseudoline arrangement on the right. }
    \label{pappas}
\end{figure}

The second problem is to count the number of actual (affine) hyperplane arrangements, but this is much harder since we don't have a perfect set-theoretical abstract model for them. We do know when the enumerations separate into two distinct problems, at least for dimensions 2 and 3.
The Pappus arrangement of 9 lines helps us find the first counterexample showing that not all oriented matroids are represented as hyperplane arrangements. This is the left-hand arrangement in Figure~\ref{pappas}, and a key feature is that Pappus' theorem, a famous theorem of  Euclidean/projective geometry, shows that the points of intersection on the top and bottom horizontal lines force the three points in the middle to be colinear. That suggests the second picture: it is an arrangement of pseudolines which  cannot be straightened while preserving their face structure---we call this situation \textit{nonstretchable}. It was found by Levi \cite{levi}, and Goodman and Pollack proved  Gr\"{u}nbaum's conjecture that no arrangement of 8 or fewer pseudolines is nonstretchable \cite{proof}. Thus the total number of  arrangements of 8 lines is 41349340. For 9 lines, the number is unknown: but Richter-Gebert showed that the example in Figure~\ref{pappas} is the unique nonstretchable case \cite{richter}.

\begin{figure} [h]
    \centering
    \includegraphics[width=.9\textwidth]{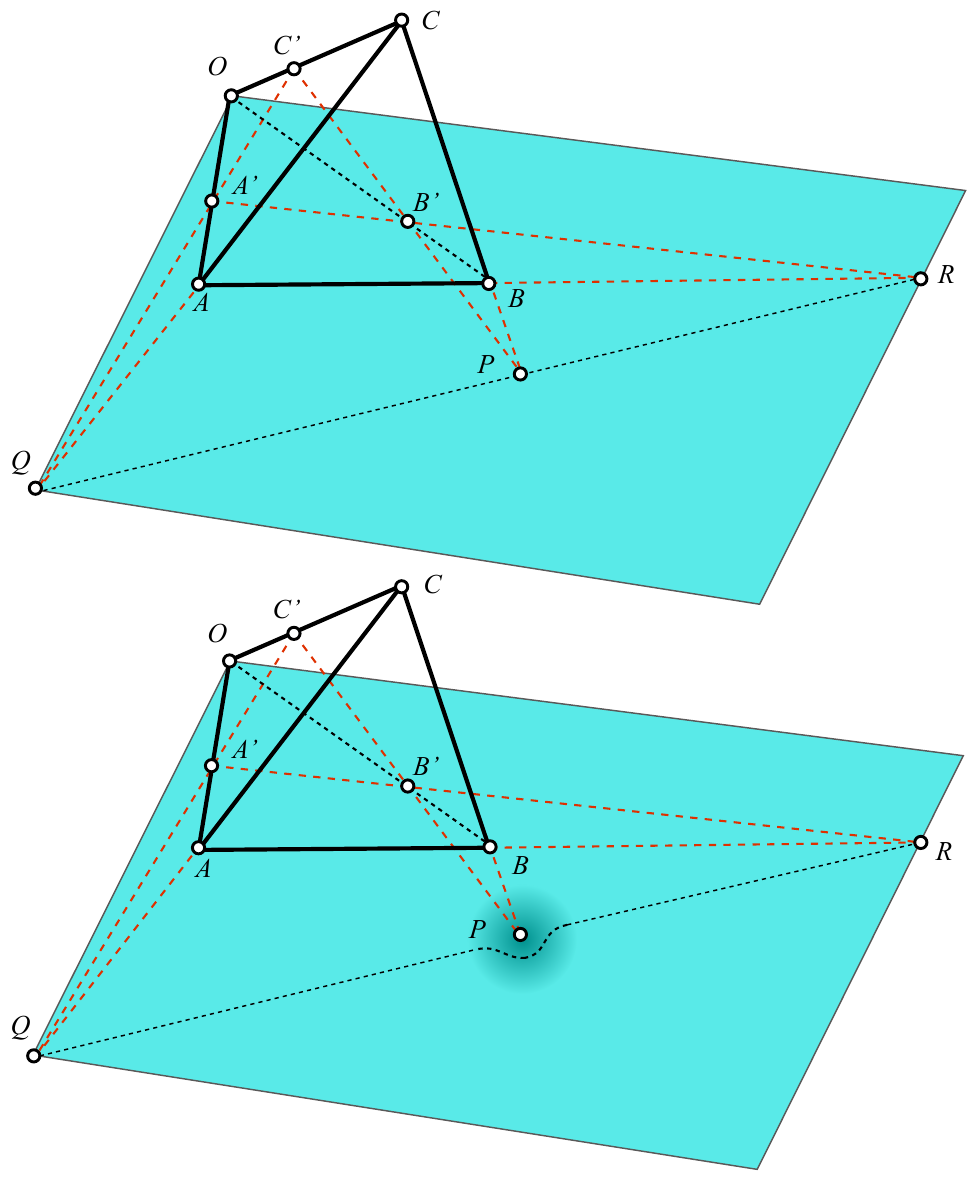}
    \caption{Top: an arrangement of 8 planes described by Goodman and Pollack in \cite{three}. Bottom: the associated nonstretchable arrangement of 8 pseudoplanes, where point $P$ is no longer on the pseudoplane containing $\Delta OQR.$}
    \label{eightnonboth}
\end{figure}

Figure~\ref{eightnonboth} shows an example of a nonstretchable pseudoplane arrangement (so an affine oriented matroid not represented by a plane arrangement) with $n=8$ pseudoplanes in dimension 3, found by Goodman and Pollack in \cite{three}.
First we pick points $A,B,C$ and $O$ in $\mathbb{R}^3$ in general position, so that they are the corners of a tetrahedron. Four of the eight planes are the four faces of the tetrahedron, extended infinitely although we only draw the triangular faces.  Then $A', B'$ and $C'$ are chosen on three edges of the tetrahedron as shown in Figure~\ref{eightnonboth}. Three more planes in our arrangement are those that contain the triangles $\Delta A'B'C$, $\Delta AB'C'$, and $\Delta A'BC'$. 
These three triangles are not drawn in Figure~\ref{eightnonboth}, but their intersections with the tetrahedral faces make the lines connecting $A', B',$ and $C'.$ Three new points are the intersections of lines:   $P$ is the intersection of lines $\overline{BC}\cap \overline{B'C'} $, $Q = \overline{AC}\cap \overline{A'C'}$, and $R = \overline{AB}\cap \overline{A'B'}$. These three points are geometrically forced to lie on a plane together with point  $O,$ (a good student exercise) and that is our eighth plane. 
Goodman and Pollack show that if we bend the eighth plane (with a dimple) to become a pseudoplane which misses the point $P$, we get a nonstretchable pseudoplane arrangement. We show that bent version in the bottom of Figure~\ref{eightnonboth}.
Goodman and Pollack also show that eight pseudoplanes are required to  make such an example, so that all arrangements of 7 or fewer planes are stretchable. Thus we can conclude by finding all the arrangements of 4 or 5 planes in $\mathbb{R}^3,$ since the total numbers coincide with the total numbers of affine oriented matroids.

\section{Illustrating plane arrangements: $n$=4 and $n$=5 planes} \label{illus}

 In general, there are several types of arrangements of $n$ planes possible in $\mathbb{R}^
3$. Here
we will describe the construction of all the affine plane arrangements in the cases where $n=4$ and  $n=5$, and illustrate some of them. Note that these constructions in general might well overlap, so that we will need to check that all the results are really distinct face-combinatorially.  Recall that central arrangements have a mutual common intersection and that essential arrangements have some subset of planes with a single point of common intersection. 

1) \textbf{The trivial arrangement} of $n$ parallel planes. This is seen for $n=4$ in the top left of Figure~\ref{nineof14}.

2) \textbf{The Cartesian product} of any nontrivial $n$-line arrangement in the $xy$-plane with the $z$-axis. There are 8 of these using 4 planes, as seen in Figure~\ref{nineof14}, using the line arrangements in Figure~\ref{fig:new8}. There are 46 of these with 5 planes, as seen via Table~\ref{counts} (where the nontrivial affine arrangements of 5 lines are the same as rank 2 affine arrangements of 5 planes). For pictures of the 46 line arrangements see Table 8.4 on page 67 of \cite{finschi}.

\begin{figure}
    \centering
    \includegraphics[width=\textwidth]{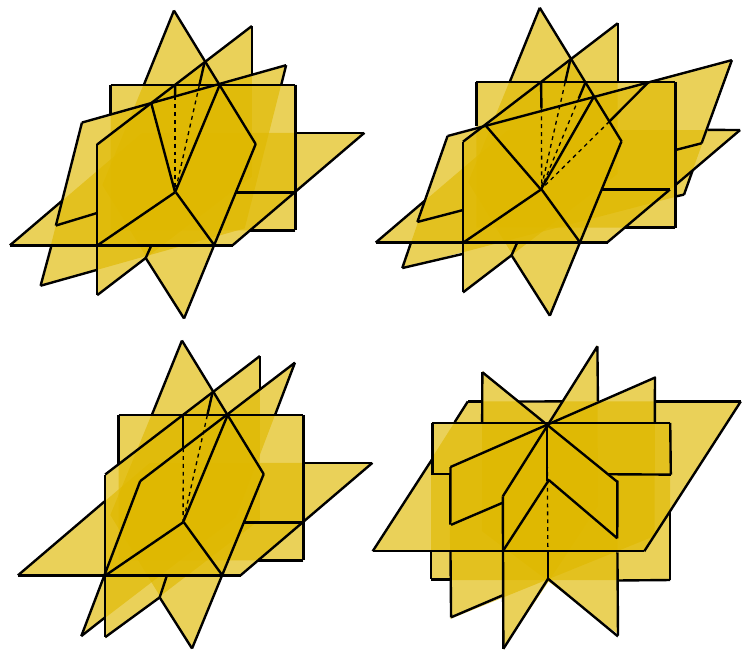}
    \caption{Four of the 27 rank 3 arrangements of five planes are central and essential, and all are found as  suspensions of 4 lines. Notice that the example in Figure~\ref{fig:fivedu} at first appears to be missing, but it is here!} 
    \label{fivesus}
\end{figure}

3) \textbf{The suspensions} of any line arrangements of $n-1$ lines.  Distinct affine line arrangements can produce equivalent  suspensions.  These are counted directly, checking for duplicates. Two of these are seen using 4 planes in Figure~\ref{fiveof14}. We find four of these using five planes, as seen in Figure~\ref{fivesus}. We leave it as an exercise to the reader to see which of the four the other line arrangements suspend to become! That there are only those four is also seen indirectly once we list the rest of the 5-plane arrangements to complete the 27 predicted by Table~\ref{counts}.
\begin{exer}
  For each of the 8 nontrivial line arrangements of 4 lines shown in Figure~\ref{fig:new8}, find the  suspension. One is shown in Figure~\ref{fig:fivedu}, but it is face-combinatorially equivalent to one of the four shown in Figure~\ref{fivesus}. Alternatively, for each of the 4 central essential plane arrangements in Figure~\ref{fivesus} find the collection of isomorphism classes of affine oriented matroids that can be formed by choosing any of the five planes to mark. Either way, find the partition of the 8 nontrivial line arrangements of 4 lines into classes determined by the feature that all in a class yield equivalent  suspensions.
\end{exer}
\begin{prob}
    Find invariants for affine hyperplane arrangements that can predict when their  suspensions will be face-combinatorially inequivalent. Test these on the 46 line arrangements of five lines shown in \cite{finschi}.
\end{prob}
4) \textbf{The bisected Cartesian product} of a non-central nontrivial $(n-1)$-line arrangement: first extended along the $z$-axis, then bisected with the $xy$-plane itself (or any plane perpendicular to the $z$-axis). There are two of these using using $n=4$ planes shown in Figure~\ref{fiveof14}. There are 7 of these using $n=5$ planes.  Figure~\ref{fiveslice} shows an example.

5) ``\textbf{New'' affine arrangements} (not using Cartesian products or suspensions).  In practice, we find most of them by focusing on bounded polytopal complexes, arrangements with various inclusions of faces that are bounded 3D chambers. 
One of these bounded chambers, the tetrahedron, is seen using 4 planes in Figure~\ref{fiveof14}. Using 5 planes  there is only one of these arrangements that has no bounded chambers, shown in Figure~\ref{fiveslice}. That leaves 15 to finish the count of 27 of rank 3 using 5 planes, all are shown in Figures~\ref{fivesingle},~\ref{fivedouble}, and~\ref{fivetriple}.

\subsection{Checking for completion} Showing these 27 rank 3 affine arrangements are all face-combinatorially distinct is straightforward, using invariants of the poset isomorphisms: any listing of the faces by dimension, boundedness, and shape. We use numbers of types of unbounded and bounded polyhedra (open and closed chambers),  numbers of triangles and quadrilaterals, numbers of types of edges (rays and line segments), and numbers of points (0-dimensional intersections). Our Figures~\ref{fivesus} --~\ref{fivetriple} are organized to have no overlap in numbers of points and bounded chambers. In Figure~\ref{fivesus} we count the numbers of rays in each: 16, 20, 12 and 10. Figure~\ref{fiveslice} includes the unique rank 3 arrangement with only one bounded edge, and an example of 7 more that can be seen to be inequivalent via Figure~\ref{fig:new8}. In Figure~\ref{fivesingle} the two arrangements with a tetrahedral chamber have 18 and 16 rays; the other two arrangements have different polyhedral chambers. In Figure~\ref{fivedouble} there are three plane arrangements with two tetrahedral chambers, but they have respectively 6, 7, and 5 points.  In that same figure there are two plane arrangements with a tetrahedron and a triangular prism, but they have 8 and 7 points respectively. Finally in that figure there is a plane arrangement with a pyramid and a tetrahedron, and another with two triangular prisms. In Figure~\ref{fivetriple} the plane arrangements with 3 chambers have different polyhedra: two tetrahedra and a prism, two tetrahedra and a square pyramid, and two prisms and a tetrahedron. For the reader we leave an easier exercise: to find invariants distinguishing all 14 of the affine arrangements of four planes in Figures~\ref{nineof14} and~\ref{fiveof14}.

\subsection{Further reading and exercises} For further elementary reading we recommend first the excellent introductions to arrangements by Halperin \cite{chaptera}, pseudolines by Goodman \cite{chapterp}, and oriented matroids by  Richter-Gebert and Ziegler \cite{chapter}, all in \cite{goodman}. More advanced reading includes  the original papers of Zaslavsky and collaborators \cite{zas, zas2, zas3, zas4} and the comprehensive monographs:  \cite{obook} and \cite{compu}. Perhaps the most exciting new developments are the specializations of oriented matroids, by adding extra requirements like \textit{purity} and \textit{positivity}, for the added value in studying networks.  The connections between pure oriented matroids and positroids are discussed by Galashin and Postnikov \cite{galp}. The connections between realizability, positively oriented matroids, and positroids are described by Ardila, Rinc\'{o}n and Williams in \cite{ard2, ard}. For more algebraic structures on the faces of hyperplanes, including some category theory and Hopf algebras, we recommend the new books by Aguiar and Mahajan: \cite{hyper} and \cite{bimon}. For  ambitious readers wanting more about positroids, check out the Amplituhedron approach to scattering matrices of Arkani-Hamed and collaborators, in \cite{ark1} and \cite{ark2}, among many other papers.

Knowing that there are precisely 74 classes of arrangements of 5 hyperplanes in $\mathbb{R}^3$  gives us power. We have drawn examples of 21; there are 6 more mentioned in Figure~\ref{fiveslice} and 47 found by extending a line arrangement along the $z$-axis (including the trivial arrangement). Now since we can show that our pictures are distinct face-combinatorially, then there are no further  possibilities.
That allows  other questions to be answered. For instance, we see that a selection of one of the 11 examples of 2 or more bounded chambers in conjunction (a complex of polyhedra) completely determines the class of a five hyperplane arrangement in $\mathbb{R}^3$.  We close with some more exercises and open problems, neither guaranteed to be easy!

\begin{figure}
    \centering
    \includegraphics[width=\textwidth]{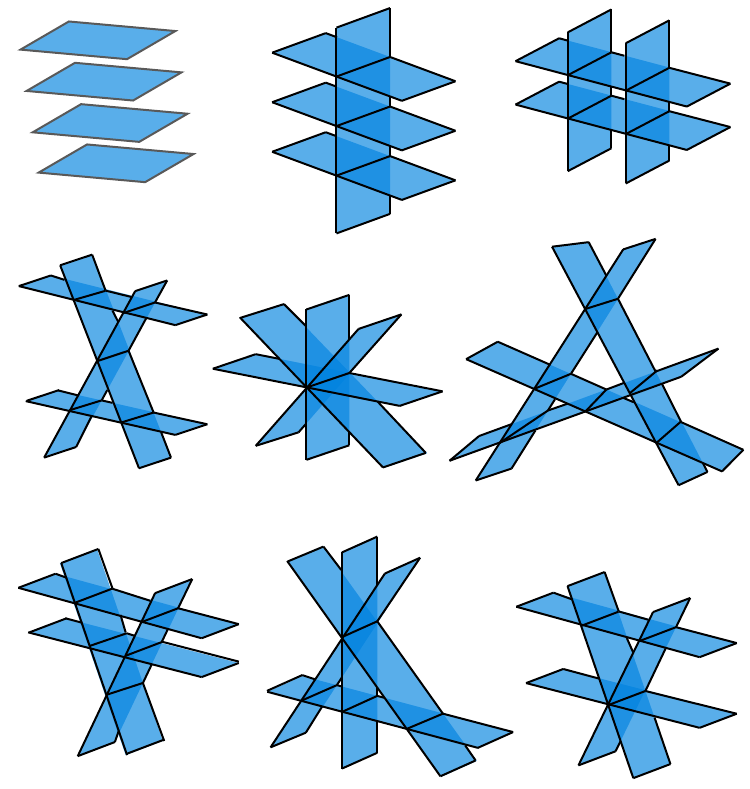}
    \caption{Nine of the 14 plane arrangements of 4 planes in $\mathbb{R}^3.$ These are the nine found as an arrangement of lines in the plane, extended along the $z$-axis. Top left is trivial (rank 1), and the rest are rank 2.}
    \label{nineof14}
\end{figure}

\begin{figure}
    \centering
    \includegraphics[width=\textwidth]{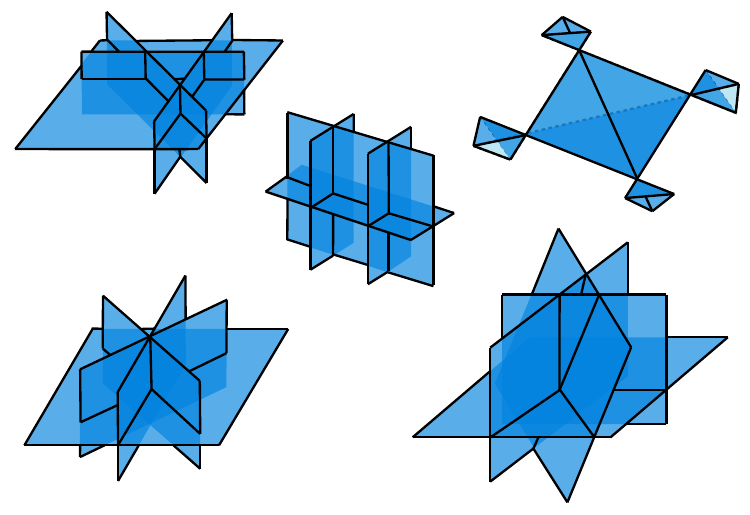}
    \caption{Five of the 14 plane arrangements of 4 planes in $\mathbb{R}^3.$ The leftmost at top and the one in the center are found by adding a perpendicular plane to a Cartesian product of a nontrivial, non-central 3-line arrangement. The rightmost on the top is a tetrahedral chamber (the four planes are extended visually a bit at each vertex). The bottommost are the two central essential arrangements, found as  suspensions of arrangements of 3 lines.}
    \label{fiveof14}
\end{figure}

\begin{exercise}
    Choosing any of the illustrated hyperplane arrangements in this paper, analyze the associated matrices. For instance, an arrangement from Figure~\ref{nineof14} or~\ref{fiveof14} comes from a matrix equation $A\mathbf{x} = \mathbf{b}$ where $A$ is $4\times 3$, $\mathbf{x}$ is a proposed solution in $\mathbb{R}^3$, and $\mathbf{b}\in\mathbb{R}^4.$ What are the rank and nullity of $A?$ Does a solution $\mathbf{x}$ exist? Can $\mathbf{b}$ be the 0-vector? Numbering the rows of $A$ for the 4 planes, which sets of rows are linearly independent? This is indeed a good way to introduce the classic \textit{matroids} as opposed to oriented matroids! Once you know what matroids are, reclassify the plane arrangements into sets of arrangements that have isomorphic matroids. 
\end{exercise}

\begin{figure}
    \centering
    \includegraphics[width=\textwidth]{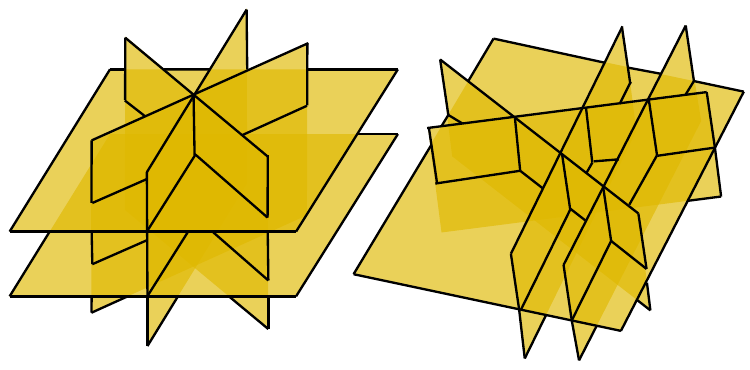}
    \caption{2 more of the 27 rank 3 arrangements of five planes. On the left is an arrangement found by adding two planes both perpendicular to a central arrangement of three planes. It is the only rank 3 non-central arrangement with exactly one bounded edge. On the right is an example of the 7 arrangements found by adding a perpendicular plane to a Cartesian product (with the $z$-axis) of one of the 7 nontrivial non-central arrangements from Figure~\ref{nineof14}. }
    \label{fiveslice}
\end{figure}

\begin{figure}
    \centering
    \includegraphics[width=\textwidth]{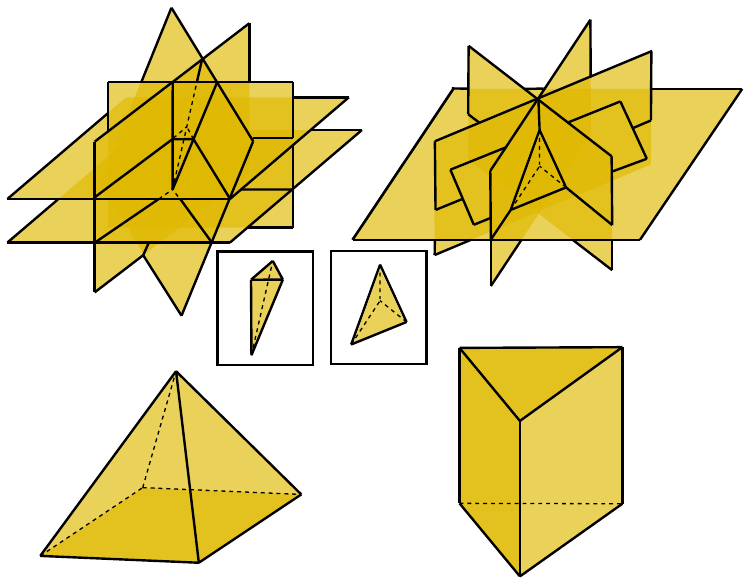}
    \caption{Four of the 27 rank 3 arrangements of five planes have a single bounded chamber: two of them have a tetrahedral chamber (shown as insets), one a square pyramid and one a triangular prism. Note that for the bottom two pictures, all five planes make bounding 2D faces (facets) of the polyhedra, so we just draw the chamber and understand the planes extend in all directions.}
    \label{fivesingle}
\end{figure}

\begin{figure}
    \centering
    \includegraphics[width=\textwidth]{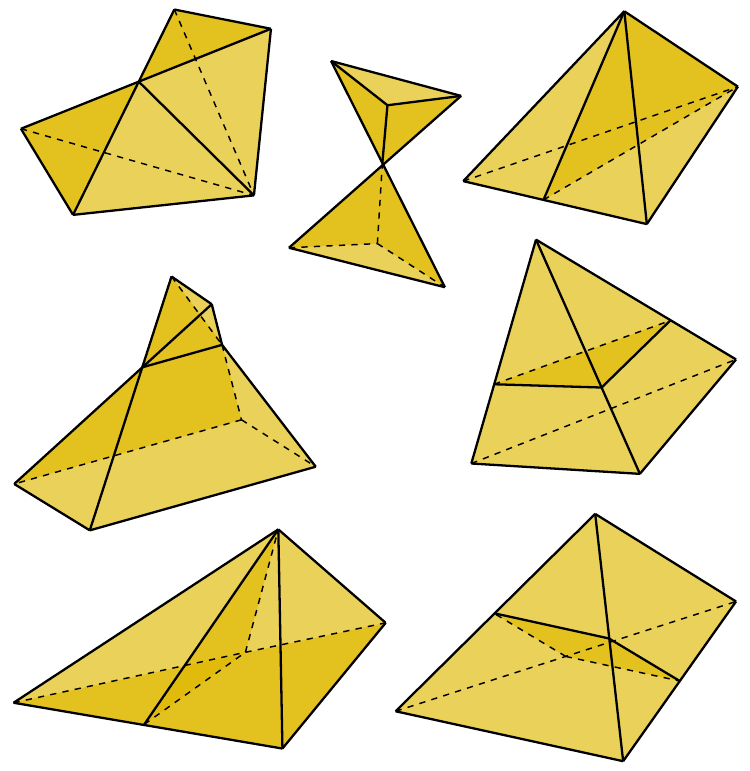}
    \caption{Seven of the 27 rank 3 arrangements of five planes have a pair of bounded chambers. Here we show only the bounded chambers of each, but exactly five planes are needed to make the bounding 2D faces of those chambers. An exercise for the reader is to find each of these by adding a single plane to one of the arrangements in Figure~\ref{nineof14} or~\ref{fiveof14}.}
    \label{fivedouble}
\end{figure}

\begin{figure}[b]
    \centering
   \includegraphics[width=\textwidth]{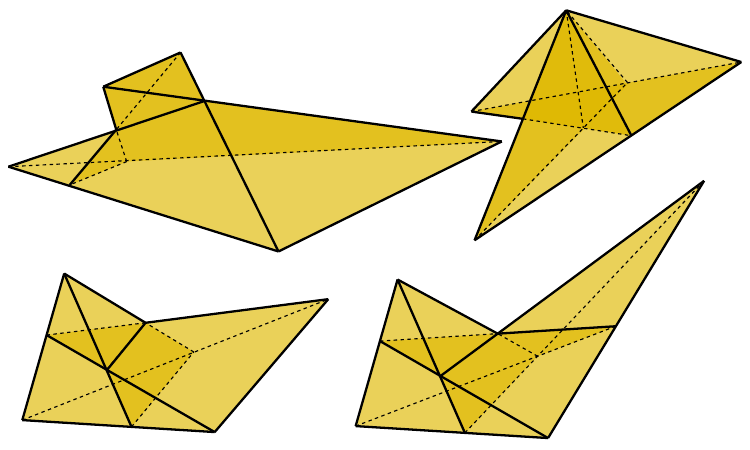}
    \caption{Three of the 27 rank 3 arrangements of five planes have a trio of bounded chambers, and the last has four bounded chambers. That last is the unique arrangement in general position: all ten subsets of three planes each determine a 0-dimensional intersection. Again we show only the bounded chambers of each, but exactly five planes are needed to make the bounding 2D faces of those chambers. An exercise for the reader is to find which of the arrangements in Figure~\ref{nineof14} or~\ref{fiveof14} can be found by deleting any single plane from one of these four arrangements.}
    \label{fivetriple}
\end{figure}

\begin{exer}
    Lots of other pictures of 5-plane arrangements can be constructed: the challenge is to always find their face-combinatorial equivalent in the listing. For instance, Figure~\ref{fiveslice} shows adding two parallel planes to one of the five arrangements in Figure~\ref{smallfive}. Do the same for any of the other four, but where is your resulting arrangement already in the list? 
\end{exer}

\begin{exer}
    The number of face-combinatorial equivalence classes of arrangements of $n$ planes that have  a single bounded tetrahedral chamber (and no other bounded chamber) is $0,1,2,\dots,$ for $n=3,4,5,\dots$. What is the number for $n=6$? Hint: notice in Figure~\ref{fivesingle} the location of the fifth plane, the one not forming a side of a tetrahedron.
\end{exer}

\begin{prob}
    The number of classes of arrangements of $n$ planes that have  a single bounded tetrahedral chamber (and no other bounded chamber) is $0,1,2,\dots,$ for $n=3,4,5,\dots$. What is the general formula for this sequence?
\end{prob}

\begin{prob}
    The number of classes of arrangements of $n$ planes that have  only a single bounded chamber (of any shape) is $0,1,4,\dots,$ for $n=3,4,5,\dots$. What is the general formula for this sequence?
\end{prob}

\begin{prob}
    The number of classes of arrangements of $n$ planes that have  only exactly two bounded chambers (of any shape) is $0,0,7,\dots,$ for $n=3,4,5,\dots$. What is the general formula for this sequence?
\end{prob}

\begin{prob}
    The number of classes of arrangements of $n$ planes that have any positive number of bounded chambers (of any shape) is $0,1,14,\dots,$ for $n=3,4,5,\dots$. What is the general formula for this sequence?
\end{prob}

We know that the maximum number of 0-dimensional faces (points) of an affine arrangement of $n$ hyperplanes in $\mathbb{R}^d$ is ${n \choose d}$. These are seen in the arrangements in general position. 
As seen in  Zaslavsky \cite{zas}, originally from Buck \cite{buck}, we have that the maximum number of chambers, both bounded and unbounded, is $\sum_{i=0}^d{n \choose i}$ (Sequence A008949 in \cite{oeis}). For instance, for $n=5$ in $\mathbb{R}^3$ we have a maximum of $1+5+10+10 = 26$ chambers of the general position plane arrangement. 


\begin{prob}
    The number of classes of  arrangements of $n$ planes that have exactly $n$ 0-dimensional faces is $0,1,4,\dots,$ for $n=3,4,5,\dots$. What is the general formula for this sequence?
\end{prob}

\begin{prob}
What is the smallest value $n$ for which parameterized equivalence of line arrangements does not give the same classes as face-combinatorial equivalence?  If $n> 8$ then the sequence A241600 can be extended by the value  4134940 as seen in \cite{finschi}.
\end{prob}

\bibliographystyle{plain}
\bibliography{plainsbib}{}

\end{document}